# THE ORIENTED SWAP PROCESS

By Omer Angel,[1] Alexander Holroyd[2] and Dan Romik

*University of Toronto, University of British Columbia and The Hebrew University*

Particles labelled $1, \ldots, n$ are initially arranged in increasing order. Subsequently, each pair of neighboring particles that is currently in increasing order swaps according to a Poisson process of rate 1. We analyze the asymptotic behavior of this process as $n \to \infty$. We prove that the space–time trajectories of individual particles converge (when suitably scaled) to a certain family of random curves with two points of non-differentiability, and that the permutation matrix at a given time converges to a certain deterministic measure with absolutely continuous and singular parts. The absorbing state (where all particles are in decreasing order) is reached at time $(2+o(1))n$. The finishing times of individual particles converge to deterministic limits, with fluctuations asymptotically governed by the Tracy–Widom distribution.

**1. Introduction.** Let $\mathcal{S}_n$ be the symmetric group of all permutations $\sigma = (\sigma(1), \ldots, \sigma(n))$ on $\{1, \ldots, n\}$, with composition denoted $(\sigma\tau)(i) := \sigma(\tau(i))$. For $1 \leq i \leq n-1$, denote the adjacent transposition or *swap* at location $i$ by $\tau_i := (i \ i+1) = (1, 2, \ldots, i+1, i, \ldots, n) \in \mathcal{S}_n$.

The *oriented swap process* is the $\mathcal{S}_n$-valued continuous-time Markov process $(\eta_t)_{t \geq 0} = (\eta_t^n)_{t \geq 0}$ defined as follows. The initial state $\eta_0$ is the *identity* permutation $\mathrm{id} := (1, 2, \ldots, n)$. From a state $\sigma$, for each $i$ satisfying $\sigma(i) < \sigma(i+1)$, the process jumps at rate 1 to the state $\sigma\tau_i$. Note that the *reverse permutation* $\mathrm{rev} := (n, \ldots, 2, 1)$ is the unique absorbing state. Our focus is on the limiting behavior of the process $(\eta_t^n)_{t \geq 0}$ as $n \to \infty$. We call the random permutation $\eta_t$ the *configuration at time $t$*, and we call $\eta_t^{-1}(k)$ the *location of particle $k$ at time $t$*. We call the function $t \mapsto \eta_t^{-1}(k)$ the *trajectory* of particle $k$; see Figures 1–3.

Received July 2008.
[1]Supported in part by an NSERC Discovery grant.
[2]Supported in part by an NSERC Discovery grant and Microsoft Research.

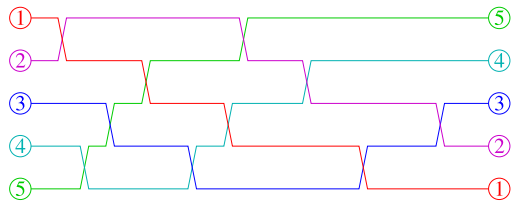

FIG. 1. *An illustration of an oriented swap process with $n = 5$. Trajectories are shown by lines.*

Our first result states that the trajectories converge to a certain family of random curves. The limiting curve for a particle at a given location is deterministic, once a random initial speed has been chosen, and is smooth, except at two points. Define the *scaled trajectory* $T_k^n = T_k : [0, \infty) \to [0, 1]$ of particle $k$ by

$$T_k^n(s) := \frac{(\eta_{ns}^n)^{-1}(k)}{n}.$$

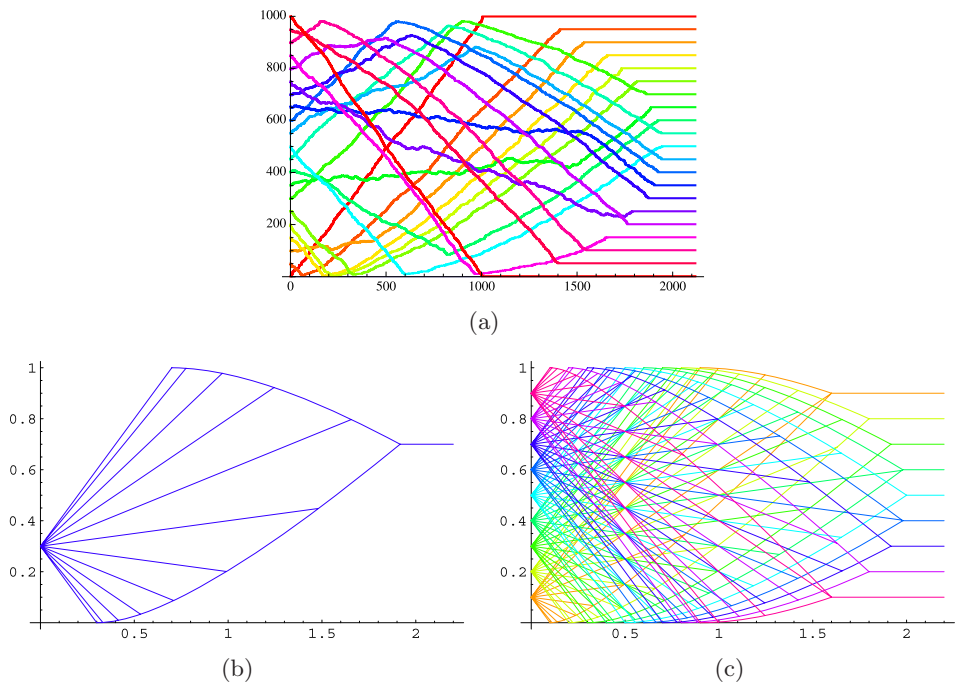

FIG. 2. (a) *Selected particle trajectories in a simulated oriented swap process with $n = 1000$;* (b) *selected possible limiting trajectories for particle $\lfloor 3n/10 \rfloor$;* (c) *selected limiting trajectories (see Theorem 1.1).*

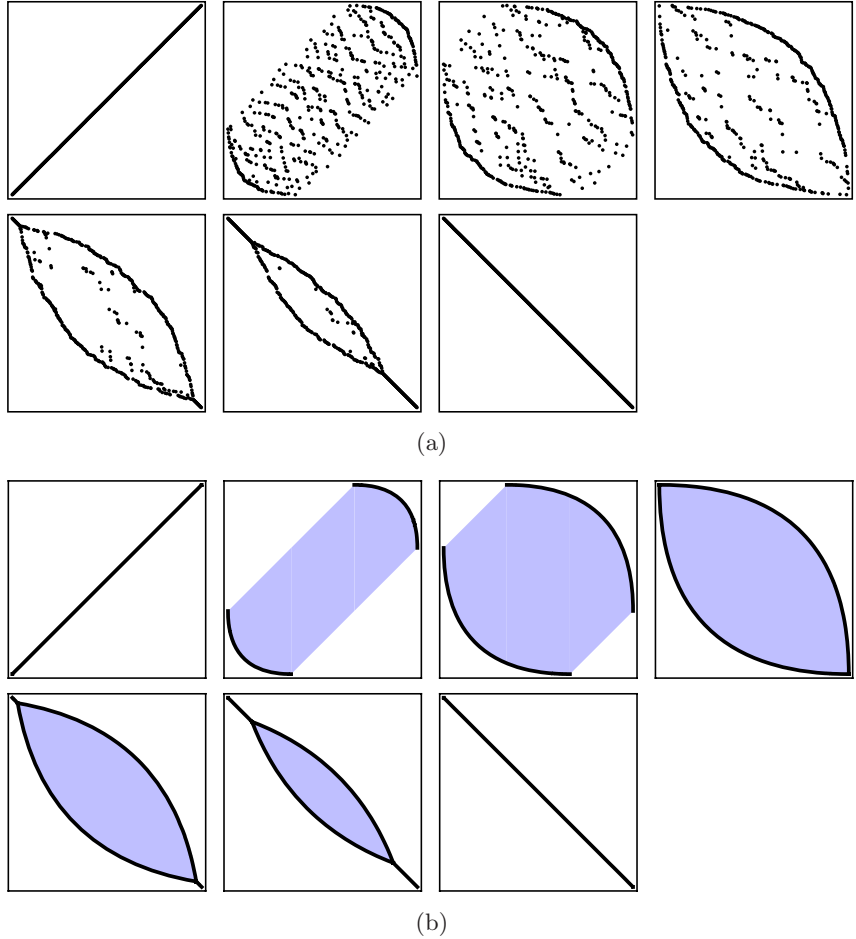

FIG. 3. (a) *The configuration (i.e., the support of the permutation matrix) at times $0, n/3, \ldots, 6n/3$ for a simulated oriented swap process with $n = 500$;* (b) *an illustration of the limiting measures for these configurations as $n \to \infty$ (see Theorem 1.2).*

THEOREM 1.1 (Trajectories). *Let $k = k(n)$ be a sequence satisfying $k/n \to y \in [0, 1]$ as $n \to \infty$. Then the scaled trajectory $T_k$ of particle $k$ satisfies*

$$T_k \Longrightarrow \phi_y \qquad as\ n \to \infty.$$

*Here, "$\Longrightarrow$" denotes convergence in distribution with respect to the uniform topology on functions $[0, \infty) \to \mathbb{R}$, and $\phi_y$ is a random function given by*

$$\phi_y(s) := \begin{cases} L_y^-(s) \vee (y + Us) \wedge L_y^+(s), & s < \gamma_y, \\ 1 - y, & s \geq \gamma_y, \end{cases}$$

where $U$ is uniformly distributed on $[-1, 1]$, and we have the deterministic functions

$$L_y^-(s) := y + s - 2\sqrt{sy}, \qquad L_y^+(s) := y - s + 2\sqrt{s(1-y)},$$

$$\gamma_y := 1 + 2\sqrt{y(1-y)}.$$

Above, $\vee$ and $\wedge$ denote maximum and minimum, respectively. Note that $a \vee b \wedge c = (a \vee b) \wedge c = a \vee (b \wedge c)$ is well-defined when $a \leq c$ and is the point in $[a, c]$ closest to $b$.

The configuration at a given time is a random permutation. The next result states that, after suitable scaling, its permutation matrix converges to a certain deterministic limit. Define the *scaled configuration* $\mu_s^n = \mu_s$ at (scaled) time $s$ to be the following random measure on $[0,1]^2$:

$$\mu_s := \frac{1}{n} \sum_{k=1}^{n} \delta\left(\frac{k}{n}, \frac{\eta_{sn}(k)}{n}\right), \tag{1}$$

where $\delta(x, y)$ denotes the point measure at $(x, y) \in \mathbb{R}^2$. The limiting measure has an absolutely continuous part, and a singular part which concentrates on a curve. For functions $g, h \colon \mathbb{R} \to \mathbb{R}$, we denote by $h(x)\delta_{y=g(x)}\, dx$ the measure on $\mathbb{R}^2$ which concentrates on the curve $\{(x,y) : y = g(x)\}$ and assigns measure $\int_A h(x)\, dx$ to $A \times \mathbb{R}$.

THEOREM 1.2 (Configurations). *For any $s > 0$, the scaled configuration $\mu_s$ satisfies*

$$\mu_s^n \Longrightarrow \kappa_s \qquad \text{as } n \to \infty.$$

*Here, "$\Longrightarrow$" denotes convergence in distribution with respect to the weak topology for random Borel measures on $\mathbb{R}^2$, and $\kappa_s$ is the deterministic measure on $[0,1]^2$ given by*

$$\begin{aligned}
\kappa_s(dx \times dy) \\
:= {} & \mathbf{1}\begin{bmatrix} L_x^-(s) < y < L_x^+(s), \\ x - s < y < x + s \end{bmatrix} \frac{1}{2s}\, dx\, dy \\
& + \mathbf{1}\begin{bmatrix} s > 1, \\ x \notin (W_s^-, W_s^+) \end{bmatrix} \delta_{y=1-x}\, dx \\
& + \left(\mathbf{1}\begin{bmatrix} s \leq 1, \\ x < s \end{bmatrix} + \mathbf{1}\begin{bmatrix} s > 1, \\ W_s^- < x < W_s^+ \end{bmatrix}\right)\left(1 - \sqrt{\frac{x}{s}}\right)\delta_{y=L_x^-(s)}\, dx \\
& + \left(\mathbf{1}\begin{bmatrix} s \leq 1, \\ x > 1 - s \end{bmatrix} + \mathbf{1}\begin{bmatrix} s > 1, \\ W_s^- < x < W_s^+ \end{bmatrix}\right)\left(1 - \sqrt{\frac{1-x}{s}}\right)\delta_{y=L_x^+(s)}\, dx,
\end{aligned} \tag{2}$$

where $L^\pm$ are as in Theorem 1.1, and for $s > 1$, we write $W_s^\pm := \frac{1 \pm \sqrt{2s - s^2}}{2}$.

The functions $W$ and $\gamma$ from Theorems 1.1 and 1.2 are inverses; more precisely, we have $\gamma_{W_s^\pm} = s$ for all $s \in [1, 2]$.

The limiting measure $\kappa_s$ in Theorem 1.2 is symmetric under the transformations $(x, y) \mapsto (1 - x, 1 - y)$ and $(x, y) \mapsto (y, x)$. The former corresponds to an obvious symmetry of the model, but the latter reflects the somewhat surprising fact that the permutation $\eta_t$ is asymptotically symmetric in law under inversion. In fact, this symmetry holds exactly for each $n$ and $t$, as the following result states.

THEOREM 1.3 (Symmetry). *For any $n$ and $t$, the permutation $\eta_t^n$ is equal in law to its inverse $(\eta_t^n)^{-1}$.*

The *inversion number* of a permutation $\sigma$ is defined as

$$\mathrm{inv}(\sigma) := \#\{(i, j) : i < j \text{ and } \sigma(i) > \sigma(j)\}.$$

An alternative description of the process $(\eta_t)$ is that the jump-rate from $\sigma$ to $\sigma\tau_i$ equals $\mathbf{1}[\mathrm{inv}(\sigma\tau_i) > \mathrm{inv}(\sigma)]$. From Theorem 1.2, we can deduce the following:

THEOREM 1.4 (Inversion number). *For each $s \geq 0$, the scaled inversion number of the configuration satisfies the convergence in probability*

$$\binom{n}{2}^{-1} \mathrm{inv}(\eta_{sn}^n) \xrightarrow{\mathbb{P}} \begin{cases} \frac{2}{3}s - \frac{1}{15}s^2, & s \in [0, 1], \\ 1 - \frac{2}{15}s^{-1/2}(2 - s)^{3/2}(2s + 1), & s \in [1, 2], \\ 1, & s \geq 2. \end{cases}$$

The limiting function on the right-hand side above is analytic, except at $s = 1$ and $s = 2$, where it is, respectively, three times and once continuously differentiable.

*Remarks on time-parameterizations.* One may consider several natural variants of the process $(\eta_t)$ in which time is parameterized differently. In the version introduced above, the total jump rate from a permutation $\sigma$ equals the size of its ascent set $A = A(\sigma) := \{i : \sigma(i) < \sigma(i + 1)\}$; we refer to this as the *variable-speed continuous-time* process. In the *fixed-speed* continuous-time process, at rate 1, we choose an $i$ uniformly from $A(\sigma)$ and jump to $\sigma\tau_i$. In the variable-speed *discrete-time* process, at each step, we choose a uniformly random $i$ from $\{1, \ldots, n - 1\}$ and jump to $\sigma\tau_i$, provided $i \in A$; while in the fixed-speed discrete-time process, we choose $i$ uniformly from $A$. Clearly, the sequence of distinct states visited has the same law for each of these four processes, and the $i$th state to be visited always has inversion number $i$. Therefore, using Theorem 1.4, one may easily translate Theorems 1.1 and 1.2 into analogous asymptotic results for the other three processes, with the limiting objects being identical, except for a deterministic time change (for brevity, we omit the full statements). Surprisingly, however, the exact symmetry in Theorem 1.3 applies only to the variable-speed models. For example, for the fixed-speed discrete-time process, with $n = 4$ and at time step 3, it is easy to check that the two mutually inverse permutations $(2, 4, 1, 3)$ and $(3, 1, 4, 2)$ have respective probabilities $1/3$ and $1/6$.

Lastly, define the *finishing time* $\beta^n(k) = \beta(k)$ of particle $k$ to be the (random) last time at which it moves, reaching its final position:

$$\beta(k) := \sup\{t > 0 : \eta_t(k) \neq n + 1 - k\}.$$

Also, define $\beta_* = \beta_*^n := \max_k \beta(k)$, that is, the hitting time of rev. The following result is strongly suggested (although not directly implied) by Theorem 1.1:

THEOREM 1.5 (Finishing times). *We have the convergence in probability*

$$\max_k |\beta^n(k)/n - \gamma_{k/n}| \xrightarrow{\mathbb{P}} 0$$

*as $n \to \infty$, where $\gamma_y := 1 + 2\sqrt{y(1-y)}$ is as in Theorem 1.1. In particular, we have*

$$\beta_*^n/n \xrightarrow{\mathbb{P}} 2.$$

Furthermore, we establish the following result on the fluctuations of the finishing times from their typical values. The *Tracy–Widom distribution function* is $F_{\mathrm{TW}}(z) := \exp[-\int_z^\infty (x-z)u(x)^2\,dx]$, where $u(x)$ is the unique solution of the Painlevé equation $u'' = 2u^3 + xu$ having the Airy function asymptotics $u(x) \sim 1/(2\sqrt{\pi})x^{-1/4}e^{-2x^{3/2}/3}$ as $x \to \infty$. The Tracy–Widom distribution originally arose [14] in random matrix theory, and has since been found to appear as a limiting law in several combinatorial models; see [3, 5, 8].

THEOREM 1.6 (Finishing time fluctuations). *Let $k = k(n)$ be such that $k/n \to y \in (0, 1)$ as $n \to \infty$. We then have the convergence in distribution*

$$\frac{\beta^n(k) - \gamma_{k/n} n}{\gamma_y^{2/3}(y(1-y))^{-1/6}n^{1/3}} \Longrightarrow F_{\mathrm{TW}} \qquad as\ n \to \infty.$$

Here, the requirement that $y \in (0, 1)$ is needed. For example, the finishing time $\beta^n(1)$ of the first particle has $\mathrm{Gamma}(n - 1, 1)$ distribution (since the particle jumps only to the right, always at rate 1), which converges, after scaling, to a Gaussian limiting law.

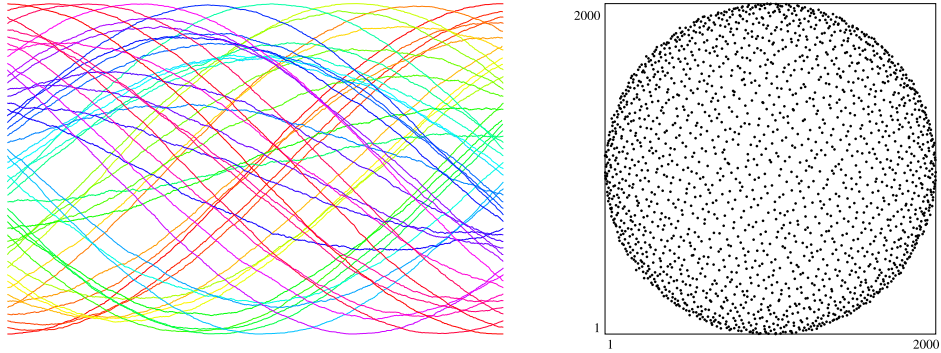

FIG. 4. *Selected trajectories and half-time configuration for a 2000-element uniform sorting network.*

*Remarks on sorting networks.* A *sorting network* is any sequence $(i_1, \ldots, i_N)$ of length $N := \binom{n}{2} = \text{inv}(\text{rev})$ such that $\tau_{i_1} \cdots \tau_{i_N} = \text{rev}$. Clearly, the sequence of swaps in the oriented swap process $(\eta_t^n)$ corresponds to a random element (with a certain non-uniform distribution) in the set of $n$-particle sorting networks. The *uniform sorting network* is instead chosen according to the uniform distribution on the same set. The present work was in part motivated by the striking results and conjectures on uniform sorting networks in [2]. The two processes behave quite differently, but share some features. For example, in the uniform sorting network, the particle trajectories conjecturally converge to random elements in a one-parameter family of curves, while the configurations conjecturally converge to a family of deterministic measures (sine curves and projected sphere measures, respectively; see Figure 4). Both of these properties hold for the oriented swap process (Theorems 1.1 and 1.2 above), but with different limiting objects.

*Remarks on the proofs.* Our analysis of the oriented swap process relies on a connection with the theory of the *totally asymmetric simple exclusion process* (*TASEP*). We show that the oriented swap process can be represented in terms of a family of coupled TASEPs, by observing that the behavior of particles $1, \ldots, k$ (if we ignore their labels) is that of a TASEP on the finite interval $[1, n]$ and then representing the TASEP on the finite interval in terms of the TASEP on $\mathbb{Z}$ using a combinatorial mapping. The main probabilistic results will follow using known limiting results for the TASEP. For Theorems 1.2 and 1.4, we use the classical hydrodynamic limit theorem of Rost [13] for the TASEP. Theorem 1.1 uses a result of Guiol and Mountford [11] on the trajectory of a second-class particle in the TASEP. Theorem 1.6 is a consequence of a theorem of Johansson [8] on the convergence of the scaled fluctuations of percolation times in oriented last-passage percolation in $\mathbb{N}^2$ to the Tracy–Widom distribution. See also [1, 6] for results related to the joint distribution of trajectories.

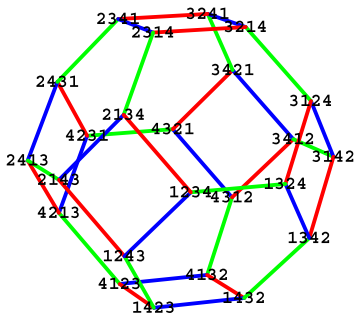

FIG. 5. *The permutahedron for $n = 4$.*

*Further remarks.* The *permutahedron* is a natural embedding in Euclidean space of the Cayley graph of $\mathcal{S}_n$ with the nearest-neighbor swaps $\tau_1, \ldots, \tau_{n-1}$ as generators. Specifically, the vertex corresponding to a permutation $\sigma$ is located at the point $\sigma^{-1} = (\sigma^{-1}(1), \ldots, \sigma^{-1}(n)) \in \mathbb{R}^n$; it is easily verified that all vertices lie on an $(n-2)$-sphere (see Figure 5). The oriented swap process may be regarded as a continuous-time simple random walk on the permutahedron, constrained to move only in directions which increase the scalar product with the vector from $\mathrm{id}^{-1}$ to $\mathrm{rev}^{-1}$.

The oriented swap process is a special case in the family of Markov processes on $\mathcal{S}_n$ in which the jump rate from $\sigma$ to $\sigma\tau_i$ is $\alpha$ if $\sigma(i) < \sigma(i+1)$, and $\alpha'$ otherwise, for fixed parameters $\alpha, \alpha'$. These processes were studied in [4], where it was proven that the mixing time for $\alpha \neq \alpha'$ is $O(n)$ [or $O(n^2)$ in the discrete-time version of the process of that paper]. Analogs of Theorems 1.2 and 1.4 for these processes can be proven along the same lines as our proofs, by using results on the *partially* asymmetric exclusion process due to Kipnis, Olla and Varadhan [9, 12] that generalize the results we used for the TASEP. Somewhat surprisingly, the symmetry condition in Theorem 1.3 also holds for these models; see [1]. However, extending Theorem 1.1 to the partially asymmetric swapping process would require proving an analogue of the Ferrari–Kipnis theorem for partially asymmetric exclusion; see Section 8.

**2. Symmetry.** In this section, we prove Theorem 1.3. As noted, the result holds for both the discrete- and continuous-time variable-speed versions of the model. The continuous-time claim follows from the discrete-time statement, which, in turn, is an immediate consequence of the following lemma.

Define the operator $S_i$, acting from the right on permutations, by

$$\sigma \cdot S_i := \begin{cases} \sigma\tau_i, & \text{if } \sigma(i) < \sigma(i+1), \\ \sigma, & \text{otherwise} \end{cases}$$

for $\sigma \in \mathcal{S}_n$; thus, $S_i$ sorts the particles in positions $i$ and $i+1$ into decreasing order. Note that we can write

$$\sigma \cdot S_i = \max\{\sigma, \sigma\tau_i\},$$

where the maximum of two permutations is the one with the greater inversion number. Note that $\operatorname{inv}(\sigma\tau_i) - \operatorname{inv}(\sigma) = \pm 1$, so the maximum above is always well-defined.

LEMMA 2.1. *For any sequence $(i_j)_{j=1}^k$, we have*

(3) $$\operatorname{id} \cdot S_{i_1} \cdots S_{i_k} = (\operatorname{id} \cdot S_{i_k} \cdots S_{i_1})^{-1}.$$

PROOF OF THEOREM 1.3. The discrete-time variable-speed process may be constructed in terms of an i.i.d. sequence of random variables $(I_j)_{j>0}$ chosen uniformly at random from $\{1, \ldots, n-1\}$. The process is then given by $\zeta_k^n = \operatorname{id} \cdot S_{I_1} \cdots S_{I_k}$. Since reversing $I_1, \ldots, I_k$ does not change their joint distribution, the theorem follows from Lemma 2.1. The continuous-time process may be defined as $\eta_t^n = \zeta_{X(t)}^n$, where $X(\cdot)$ is the counting function for a Poisson process with rate $n-1$ which is independent of the swap locations. The result follows. □

PROOF OF LEMMA 2.1. We first prove the following claim. Let $i, i_1, \ldots, i_n$ be any integer sequence. If we let $X = \operatorname{id} \cdot S_{i_1} \cdots S_{i_n}$ and $Y = \tau_i \cdot S_{i_1} \cdots S_{i_n}$, then

$$Y^{-1} = \max\{X^{-1}, X^{-1}\tau_i\} = X^{-1} \cdot S_i.$$

Indeed, it is easy to see, by induction on $k$, that for any $k \leq n$ and any particle $j \notin \{i, i+1\}$, the locations of $j$ in $\operatorname{id} \cdot S_{i_1} \cdots S_{i_k}$ and in $\tau_i \cdot S_{i_1} \cdots S_{i_k}$ are equal, since particle $j$ does not distinguish between $i$ and $i+1$ when attempting a swap. For $k = n$, this shows that $Y^{-1}$ is one of $X^{-1}$ and $\tau_i X^{-1}$. However, in $Y$, we must have $Y^{-1}(i+1) < Y^{-1}(i)$, since that is so in $\tau_i$ and particles $i$ and $i+1$ cannot be swapped again when applying $S_{i_1}, \ldots, S_{i_n}$.

The lemma is now proved by induction. Assume it is true for some sequence $i_1, \ldots, i_n$ and extend the sequence by adding $i_0 = i$. We then have, applying the induction hypothesis for $X$ (with $X, Y$ as above), that

$$(\operatorname{id} \cdot S_{i_0} \cdots S_{i_n})^{-1} = Y^{-1} = X^{-1} \cdot S_{i_0} = \operatorname{id} \cdot S_{i_n} \cdots S_{i_0},$$

as claimed. □

**3. The infinite oriented swap process and exclusion processes.** While our results are concerned with processes taking place within a finite interval $[1, n]$, in order to prove them, it will be useful to see the interval as part of $\mathbb{Z}$ and the finite-interval process as a function of a process on the entire line. We now introduce this process, which will be called the *infinite oriented swap process*.

In the infinite oriented swap process, particles with labels from $\mathbb{Z}$ occupy positions in $\mathbb{Z}$, with each position containing exactly one particle at any given time, so the configuration space is $\mathbb{Z}^{\mathbb{Z}}$. The infinite oriented swap process is the $\mathbb{Z}^{\mathbb{Z}}$-valued continuous-time Markov process $(\zeta_t)_{t \geq 0}$ defined as follows. The initial state $\zeta_0$ is the infinite identity configuration id [defined by $\mathrm{id}(k) = k$ for all $k \in \mathbb{Z}$]. For each edge $(k, k+1)$, the particles in positions $k$ and $k+1$ "attempt to swap" at rate 1, succeeding if and only if they are in increasing order.

The existence of this process is proved using a graphical representation. Specifically, for each $k \in \mathbb{Z}$, let $\Pi_k$ be a Poisson process with density 1 on $\mathbb{R}^+ = [0, \infty)$, where $(\Pi_k)_{k \in \mathbb{Z}}$ form an independent family. For each $k \in \mathbb{Z}$, let $\Pi_k$ be the set of times at which a swap is attempted between positions $k$ and $k+1$. Since almost surely, for all $k \in \mathbb{Z}$ and $t > 0$, we have that $|\Pi_k \cap [0, t]| < \infty$, and, a.s. for all $t > 0$, we have $|\Pi_j \cap [0, t]| = 0$ for some arbitrarily large (both positive and negative) values of $j$, it follows that the label of the particle in position $k$ and the location of particle $j$ at any time $t$ are well-defined.

As before, we call $\zeta_t$ the *configuration at time $t$*, we call $\zeta_t^{-1}(k)$ the *location of particle $k$* at time $t$ and we call the function $t \mapsto \zeta_t^{-1}(k)$ the *trajectory* of particle $k$.

Note that for finite $n$, the oriented swap process $(\eta_t^n)_{t \geq 0}$ can be realized similarly to the infinite process by using only the Poisson processes $\{\Pi_k\}_{1 \leq k \leq n-1}$ and ignoring all the others. More generally, for any (possibly infinite) interval $I = [a, b] \subseteq \mathbb{Z}$, where $a < b$ and $a, b \in \mathbb{Z} \cup \{\pm \infty\}$, we consider the oriented swap process on $I$ defined by restricting configurations to $I$ and applying only swaps coming from the Poisson processes $\{\Pi_k\}_{k \in [a, b-1]}$. We denote this process by $(\zeta_t^I)_{t \geq 0}$, with the conventions that if $I = \mathbb{Z}$, then the superscript $I$ is omitted, and that $\zeta_t^n = \zeta_t^{[1,n]} = \eta_t^n$.

The *totally asymmetric simple exclusion process* (*TASEP*) is a process on $\mathbb{Z}$ (or, for our purposes in some cases, a subinterval of $\mathbb{Z}$) with just two kinds of particles (called a "particle" and a "hole"), where, with exponential rate 1, each particle tries to jump one step to the right (leaving a hole behind), succeeding if the place to its right contains a hole. On a finite interval, particles cannot leave the interval from the right or enter it from the left, and the process will support only a finite number of moves before the particles are stuck at the right side of the interval.

It is well-known and easy to see that a TASEP on $I = [a, b]$ can be constructed starting from an arbitrary initial configuration using a family of independent Poisson processes $\{\Pi_k\}_{k \in [a, b-1]}$ of attempted jump times, similarly to the construction of the oriented swap processes above. We shall be interested in TASEPs with a particular class of step functions as the initial conditions. For any interval $I = [a, b] \subseteq \mathbb{Z}$ and $k \in \mathbb{Z}$, let $(\nu_t^{k,I})_{t \geq 0}$ denote the TASEP on the interval $I$ with initial condition

$$\nu_0^{k,I}(x) = \mathbf{1}_{\{x \leq k\}}$$

constructed from the same infinite family of Poisson processes $\{\Pi_k\}_{k \in \mathbb{Z}}$ that was used in the construction of the infinite oriented swap process above. As before, if $I = \mathbb{Z}$, we may omit $I$ from the superscript and denote the process simply as $\nu_t^k$, and if $I = [1, n]$, we denote the corresponding process by $\nu_t^{k,n}$.

For any $k \in \mathbb{Z}$ and a configuration $\rho \in \mathbb{Z}^I$, we define $T_k \rho \in \{0, 1\}^I$ by

$$(T_k \rho)(x) = \mathbf{1}_{\{\rho(x) \leq k\}},$$

that is, the composition of the characteristic function of $(-\infty, k]$ with $\rho$.

LEMMA 3.1. *With the above construction of the processes, a.s. for any subinterval $I \subseteq \mathbb{Z}$ and $t \geq 0$, we have*

$$\nu_t^{k,I} = T_k \zeta_t^I.$$

PROOF. The identity holds initially by definition and is preserved by any attempted swap. $\square$

Our next goal is to describe the relations between $(\nu_t^{k,I})_{t \geq 0}$ and $(\nu_t^k)_{t \geq 0} = (\nu^{k,\mathbb{Z}})_{t \geq 0}$. Informally, if $I \subset J \subseteq \mathbb{Z}$ are two subintervals of $\mathbb{Z}$ (possibly infinite), then the randomness involved in $(\nu_t^{k,I})_{t \geq 0}$ is a subset of the randomness involved in $(\nu_t^{k,J})_{t \geq 0}$. It turns out that $\nu_t^{k,I}$ is a function of $\nu_t^{k,J}$, not just as a process, but also at any fixed time $t$.

For what follows, it will be more convenient to work exclusively with configurations in $\{0, 1\}^{\mathbb{Z}}$ having a rightmost particle. Let $\Omega_0 \subset \{0, 1\}^{\mathbb{Z}}$ be the space of such configurations. If $\rho \in \{0, 1\}^I$, extend $\rho$ to $\{0, 1\}^{\mathbb{Z}}$ by setting $\rho(x) = 0$ for $x \in \mathbb{Z} \setminus I$. Note that $\Omega_0$ is a.s. invariant under the Markovian dynamics of the TASEP, so, because of our choice of initial conditions, we get that a.s. $\nu_t^{k,I} \in \Omega_0$ for all $t > 0, k \in \mathbb{Z}$ and $I \subset \mathbb{Z}$.

Define the operators $R_k$, $B_n$ and $J_m$ (where $m \in \mathbb{Z}$) on $\Omega_0$ as follows. We call $R_k$ the *cut-off operator*. For a configuration $\rho$, $R_k(\rho)$ keeps only the $k$ rightmost particles of $\rho$. Formally,

$$(R_k \rho)(x) = \begin{cases} \rho(x), & \text{if } \sum_{y > x} \rho(y) < k, \\ 0, & \text{otherwise.} \end{cases}$$

The *push-back operator* $B_n$ pushes all particles back into $(-\infty, n]$ and preserves the exclusion. $B_n$ moves the $j$th rightmost particle, if it is in location $x$, to location $x \wedge (n+1-j)$. Formally,

$$(B_n\rho)(x) = \begin{cases} 0, & x > n, \\ 1, & \text{if } x \leq n \text{ and } \sum_{y>x}\rho(y) > n-x, \\ \rho(x), & \text{if } x \leq n \text{ and } \sum_{y>x}\rho(y) \leq n-x. \end{cases}$$

The *jump operator* $J_m$ is analogous to the sorting operator $S_i$ from the previous section and tries to make a particle at $m$ jump to $m+1$ if there is a particle at $m$ and no particle at $m+1$. Formally,

$$J_m\rho = \begin{cases} \rho \cdot \tau_m, & \text{if } \rho(m) = 1 \text{ and } \rho(m+1) = 0, \\ \rho, & \text{otherwise,} \end{cases}$$

where $\rho \cdot \tau_m$ denotes $\rho$ with the values at $m$ and $m+1$ transposed.

One more notion that will prove useful is the *queue-length function*. For any $x \in \mathbb{Z}$ and a configuration $\rho \in \Omega_0$, $S(\rho, x)$ gives the number of particles to the right of $x$:

$$S(\rho, x) = \sum_{y>x} \rho(y).$$

Clearly, a configuration $\rho \in \Omega_0$ is completely determined by $S(\rho, \cdot)$. It is easy to see that in terms of the $S$-function, the operators $R_k$ and $B_n$ take the following form:

$$S(R_k\rho, x) = S(\rho, x) \wedge k, \tag{4}$$

$$S(B_n\rho, x) = S(\rho, x) \wedge (n-x)^+ \tag{5}$$

(where, for $u \in \mathbb{R}$, we use the notation $u^+ = u \vee 0$). It follows that $B_n$ and $R_k$ commute. We will also need the following result:

LEMMA 3.2.  (i) *If* $1 \leq m \leq n-1$, *then*

$$B_n R_k J_m = J_m B_n R_k.$$

(ii) *If* $m \geq n$, *then*

$$B_n R_k J_m = B_n R_k.$$

(iii) *If* $m \leq 0$ *and if* $\rho \in \Omega_0$ *has its $k$ rightmost particles in* $[1, \infty)$, *then*

$$B_n R_k J_m \rho = B_n R_k \rho.$$

PROOF. We first prove (ii). $B_n$ and $R_k$ commute, so it suffices to show that $B_n J_m = B_n$ whenever $m \geq n$. This is true since $B_n \rho$ depends only on $S(\rho, k)$ for $k < n$, whereas $J_m$ does not change any of those numbers.

Claim (iii) is similarly easy. $R_k \rho$ depends only on the positions of the $k$ rightmost particles in $\rho$. Thus, if there are at least $k$ particles in $[1, \infty)$ and $m \leq 0$, then $J_m$ does not affect any of them and $R_k J_m \rho = R_k \rho$. Claim (iii) follows.

To prove (i), we show that for $m \in [1, n-1]$, $J_m$ commutes with both $B_n$ and $R_k$. First, we prove that $R_k J_m = J_m R_k$ for any $m \in \mathbb{Z}$. Imagine that the $k$ rightmost particles in a configuration $\rho \in \Omega_0$ are colored red and that all other particles are colored blue. Then the statement that $R_k J_m \rho = J_m R_k \rho$ simply says that making a particle at $m$, if there is one there, try to jump (note that its color will be preserved whether it jumps or not) and then deleting all blue particles is the same as first deleting all blue particles and then trying a jump at $m$. If $\rho$ has no particle at $m$, then this is clearly true, and if it has a particle, then the statement can be false only if the particle is a red particle which has a blue particle to its right. But this cannot happen since all blue particles are to the left of all red particles.

Next, we prove that if $m < n$, then $B_n J_m = J_m B_n$. If, in a configuration $\rho$, there is no particle at $m$, then $J_m$ leaves $\rho$ unchanged. Then, if, in $B_n \rho$, there is also no particle at $m$, we are done; otherwise, there is a particle that was pushed there by other particles from the right, so, in particular, $B_n \rho$ also has a particle at $m+1$ (here, we use the fact that $m < n$), and $J_m$ also leaves $B_n \rho$ unchanged, so we are done.

Alternatively, assume that $\rho$ has a particle at $m$. Let $j$ denote the ranking of that particle (in terms of right-to-left order of appearance) and let $m'$ denote the location of the $j$th rightmost particle in $B_n \rho$ (we think of it as the "same" particle after the push-back operation). If $m' < m$, then the particle was pushed, so $m' = n + 1 - j$. Applying $J_m$ will leave $B_n \rho$ unchanged since, in $B_n \rho$, there are particles at all locations between $m'$ and $n$, including $m+1$. In the other direction, applying $J_m$ to $\rho$ might make the $j$th particle jump to $m+1$, but, after applying $B_n$, it will again be pushed to $m' = n + 1 - j$. So, in this case, too, we have shown that $B_n J_m \rho = B_n \rho = J_m B_n \rho$ (note that all other particles are unaffected by $J_m$). It remains to deal with the case $m' = m$: in this case, if $m = n + 1 - j$, then whether applying $J_m$ to $\rho$ produces a jump or not, the $j$th particle will be pushed back by $B_n$ to $n + 1 - j$, with a particle in front of it blocking a jump (since $m < n$). And, if $m < n + 1 - j$, then the $j$th particle does not get pushed back, so, if $\rho$ has a particle at $m+1$, that particle will still be there after applying $B_n$ and, in both cases, a jump will not occur; otherwise, if there is no particle at $m+1$, after applying $B_n$, the $(j-1)$th rightmost particle will be pushed to $n + 2 - j > m + 1$, so, after applying $B_n$, there will still not be a particle at $m+1$ and a jump will occur for both orders of applying the operations. This completes the proof that $B_n J_m = J_m B_n$ and hence the proof of (i). □

LEMMA 3.3. *Let $k, n \in \mathbb{Z}$. Almost surely, for any $t \geq 0$, we have*

$$\nu_t^{k,n} = B_n R_k \nu_t^k. \tag{6}$$

PROOF. The equality (6) is satisfied at time $t = 0$. Attempted jumps outside of $[1, n-1]$ have no effect on $\nu_t^{k,n}$ and, by parts (ii) and (iii) of Lemma 3.2, these have no effect on $B_n R_k \nu_t^k$ either. Jumps inside $[1, n-1]$ a.s. occur at a discrete set of times and, by Lemma 3.2(i), whenever a jump is attempted between positions $(m, m+1)$, if (6) was satisfied before the attempted jump, it will remain true after it. □

**4. Hydrodynamic limits.** The TASEP $(\nu_t^0)_{t \geq 0}$ with initial condition $\nu_0^0 = \mathbf{1}_{(-\infty,0]}$ has been studied in great depth and is the simplest case of a *shock* in a TASEP. Lemma 3.3 allows us to tap into this knowledge. We use the following fundamental result of Rost ([13], Theorem 1; see also [10], Chapter VIII, Section 5):

THEOREM 4.1 (Rost [13]). *For any $-\infty \leq u < v \leq \infty$, we have that a.s.*

$$\lim_{t \to \infty} \frac{1}{t} \sum_{ut < j < vt} \nu_t^0(j) = \int_u^v h(x) \, dx,$$

*where*

$$h(x) = 1 \wedge \frac{1-x}{2} \vee 0 = \begin{cases} 1, & x < -1, \\ \frac{1-x}{2}, & -1 \leq x \leq 1, \\ 0, & x > 1. \end{cases}$$

The following immediate corollary is an equivalent formulation of Rost's result which can be interpreted as saying that the family of functions $(x \to h(x/s))_{s \geq 0}$ is the limiting time-evolution, or *hydrodynamic limit*, of the density profile of the process $(\nu_t^0)_{t \geq 0}$, when one scales both the time- and space-axes by a parameter $n$ that goes to infinity.

COROLLARY 4.2. *For any $-\infty \leq u < v \leq \infty$ and $s > 0$, we have a.s.*

$$\lim_{n \to \infty} \frac{1}{n} \sum_{nu < j < nv} \nu_{ns}^0(j) = \int_u^v h\left(\frac{x}{s}\right) dx.$$

With this preparation, we can now formulate and prove an analogous hydrodynamic limit theorem for the TASEP $\nu_t^{k,n}$, where $k \approx y \cdot n$ for some fixed $0 < y < 1$. In terms of the original oriented swap process $(\eta^n_t)_t$, this describes the limiting flow of particles with label $\leq k$, which start out in positions $[1, k]$, to their final positions at $[n+1-k, n]$, as a function of the scaled time parameter. This result is essentially an encoded form of Theorem 1.2 and we will derive Theorem 1.2 from it later.

THEOREM 4.3. *Fix $0 < y < 1$, $s > 0$ and $0 \leq u < v \leq 1$. Let $k = k(n)$ be a sequence of integers such that $k/n \to y$ as $n \to \infty$. Then, the number of particles in the oriented swap process $(\eta^n_t)_t$ with index $\leq k$ that are in the interval $(nu, nv)$ at time $ns$ satisfies*

$$(7) \quad \frac{1}{n}\#\{1 \leq j \leq k : (\eta^n_{ns})^{-1}(j) \in (nu, nv)\} \xrightarrow[n \to \infty]{\mathbb{P}} \int_u^v f(s, x, y)\, dx,$$

*where $f$ is the function of $s, x$ and $y$ defined by*

$$(8) \quad f(s, x, y) = \begin{cases} \dfrac{s + y - x}{2s}, & (y - s) \vee L_y^-(s) < x < (y + s) \wedge L_y^+(s), \\ 1, & 0 < x < y - s, \\ 0, & y + s < x < 1, \\ 1, & s > 1 - y \text{ and } (1 - y) \vee L_y^+(s) < x < 1, \\ 0, & s > y \text{ and } 0 < x < (1 - y) \wedge L_y^-(s), \end{cases}$$

*where $L_y^{\pm}(s)$ are as defined in Theorem 1.1 in the Introduction.*

Figure 6 shows two convenient ways of visualizing the limiting density profile $f(s, x, y)$ for fixed $y$ (in this case, $y = 0.3$): in Figure 6(a), we see $f(s, x, y)$ in the $(s, x)$-plane; Figure 6(b) shows a succession of plots of $f(s, x, y)$ as a function of $x$, for several increasing values of $s$. If we think of particles with label $\leq k$ in the oriented swap process as "red" particles and particles with label $> k$ as "green," this illustrates how the red particles advance into the green zone, eventually displacing all green particles in positions $[n+1-k, n]$.

PROOF OF THEOREM 4.3. It is clearly enough to prove (7) for $v = \infty$. If $k/n \to y$ as $n \to \infty$, then, by the definitions,

$$\frac{1}{n}\#\{1 \leq j \leq k : (\eta^n_{ns})^{-1}(j) \in (nu, \infty)\} = \frac{1}{n} S(\nu^{k,n}_{ns}, nu).$$

By Lemma 3.3, this is a.s. equal to

$$\frac{1}{n}[S(\nu^k_{ns}, nu) \wedge k \wedge (n - \lfloor nu \rfloor)^+].$$

By Corollary 4.2, together with translation-invariance of the TASEP dynamics, as $n \to \infty$, this last quantity converges in probability to

$$(9) \quad F(s, u, y) := \int_u^\infty h\left(\frac{x - y}{s}\right) dx \wedge y \wedge (1 - u)^+.$$

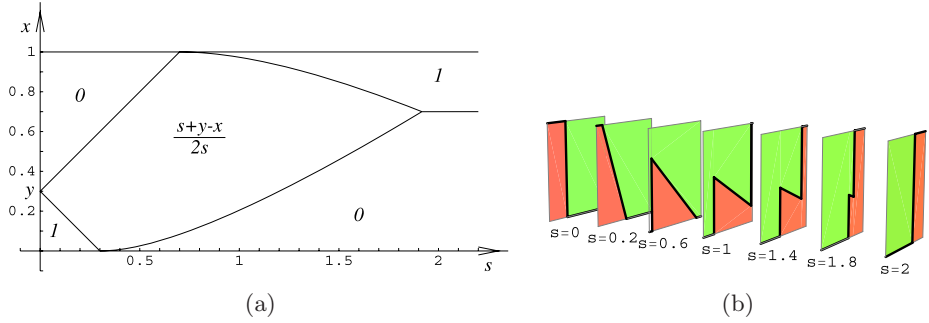

FIG. 6. (a) *Schematic representation of $f(s,x,y)$ in the $(s,x)$-plane for $y=0.3$; (b) several time slices of $f(s,x,y)$ (shown as a function of $x$) for $y=0.3$.*

It now remains to show that

$$F(s,u,y) = \int_u^\infty f(s,x,y)\,dx, \tag{10}$$

where $f(s,x,y)$ is defined in (8). Define the functions

$$\Lambda_y^-(s) = \inf\left\{u : \int_u^\infty h\left(\frac{x-y}{s}\right) dx \geq y\right\},$$

$$\Lambda_y^+(s) = \inf\left\{u : \int_u^\infty h\left(\frac{x-y}{s}\right) dx \geq (1-u)^+\right\}.$$

We can evaluate $\Lambda_y^-$ and $\Lambda_y^+$ explicitly, as shown in Figure 7, to give

$$\Lambda_y^-(s) = \begin{cases} 0, & 0 \leq s < y, \\ L_y^-(s), & s \geq y, \end{cases}$$

$$\Lambda_y^+(s) = \begin{cases} 1, & 0 \leq s < 1-y, \\ L_y^+(s), & s \geq 1-y \end{cases}$$

[e.g., from Figure 7, it is easy to see that $\Lambda_y^-(s)$ for $s \geq y$ is the solution of the equation $(y+s-\Lambda)^2/2s = y$; we omit the detailed verification of the above formulae].

Now, check that $\Lambda_y^-(s) \leq \Lambda_y^+(s)$ if $0 \leq s \leq \gamma_y$ and $\Lambda_y^-(s) > \Lambda_y^+(s)$ if $s > \gamma_y$. Equipped with this information, it is easy to write the following explicit formulas for $F(s,u,y)$. Assuming that $y \leq 1/2$, we get

$$F(s,u,y) = \begin{cases} y-u, & 0 \leq u \leq y-s, \\ \dfrac{(s+y-u)^2}{4s}, & y-s \leq u \leq y+s, \\ 0, & y+s \leq u \leq 1 \end{cases} \quad (\text{if } 0 \leq s \leq y),$$

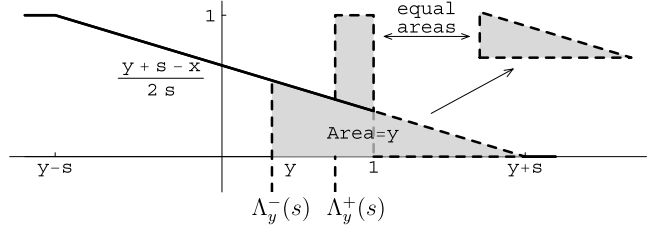

FIG. 7. *The functions* $\Lambda_y^-(s), \Lambda_y^+(s)$.

$$F(s,u,y) = \begin{cases} y, & 0 \leq u \leq L_y^-(s), \\ \dfrac{(s+y-u)^2}{4s}, & L_y^-(s) \leq u \leq y+s, \\ 0, & y+s \leq u \leq 1 \end{cases} \quad \text{(if } y < s \leq 1-y\text{)},$$

$$F(s,u,y) = \begin{cases} y, & 0 \leq u \leq L_y^-(s), \\ \dfrac{(s+y-u)^2}{4s}, & L_y^-(s) \leq u \leq L_y^-(s), \\ 1-u, & L_y^-(s) \leq u \leq 1 \end{cases} \quad \text{(if } 1-y < s \leq \gamma_y\text{)},$$

$$F(s,u,y) = \begin{cases} y, & 0 \leq u \leq 1-y, \\ 1-u, & 1-y \leq u \leq 1 \end{cases} \quad \text{(if } s > \gamma_y\text{)}.$$

On the other hand, referring to Figure 6 for convenience, one may rewrite $f(s,x,y)$ defined in (8) more explicitly (again in the case $y \leq 1/2$) as

$$f(s,x,y) = \begin{cases} 1, & 0 < x < y-s, \\ \dfrac{s+y-x}{2s}, & y-s < x < y+s, \\ 0, & y+s < x < 1 \end{cases} \quad \text{(if } 0 \leq s \leq y\text{)},$$

$$f(s,x,y) = \begin{cases} 0, & 0 < x < L_y^-(s), \\ \dfrac{s+y-x}{2s}, & L_y^-(s) < x < y+s, \\ 0, & y+s < x < 1 \end{cases} \quad \text{(if } y < s \leq 1-y\text{)},$$

$$f(s,x,y) = \begin{cases} 0, & 0 < x < L_y^-(s), \\ \dfrac{s+y-x}{2s}, & L_y^-(s) < x < L_y^-(s), \\ 1, & L_y^-(s) < x < 1 \end{cases} \quad \text{(if } 1-y < s \leq \gamma_y\text{)},$$

$$f(s,x,y) = \begin{cases} 0, & 0 < x < 1-y, \\ 1, & 1-y < x < 1 \end{cases} \quad \text{(if } s > \gamma_y\text{)}.$$

Comparing the two sets of formulae, it is clear that (10) holds. The case $1/2 < y \leq 1$ is dealt with similarly and is omitted. □

PROOF OF THEOREM 1.2. To prove (2), it is enough to prove that for each $s > 0$ and $0 \leq x, y \leq 1$, we have that
$$\mu_s^n([0,x] \times [0,y]) \xrightarrow[n\to\infty]{\mathbb{P}} \kappa_s([0,x] \times [0,y]).$$
However, looking at the definitions and using Theorem 4.3, we see that
$$\mu_s^n([0,x] \times [0,y]) = \frac{1}{n}\#\{1 \leq j \leq ny : (\eta_{ns}^n)^{-1}(j) \in [0, nx]\}$$
$$\xrightarrow[n\to\infty]{\mathbb{P}} \int_0^x f(s, u, y)\, du.$$
It remains to verify that
$$\int_0^x f(s,u,y)\, du = \kappa_s([0,x] \times [0,y])$$
or, equivalently, that
$$(11) \qquad f(s,x,y) = \frac{\partial}{\partial x} \kappa_s([0,x] \times [0,y]) =: g(s,x,y).$$
Note that, for $s \geq 2$, $\kappa_s$ is simply the arc-length measure on the segment $\{(x,y) \in [0,1]^2 : y = 1-x\}$, normalized to be a probability measure, so the identity is trivial to check. There are two other cases, $0 < s \leq 1$ and $1 < s < 2$, and in each of these cases, it is easy to compute $g(s,x,y)$ by dividing the unit square into the various possibilities for $(x,y)$. Figures 8(a) and (b) show the result (see also Section 6, where we give an alternative geometric description of $\kappa_s$). Comparing this to the explicit formulas for $f(s,x,y)$ in the case $y \leq 1/2$, one verifies that (11) holds in this case. The case $1/2 < y \leq 1$ may be dealt with similarly or can be deduced from the symmetry of the process with respect to reversing left and right, and replacing particle $k$ by $n+1-k$. □

## 5. Finishing times.

PROOF OF THEOREM 1.5. To analyze the time $\beta^n(k)$ at which particle $k$ completes its movement, note that, in the oriented swap process $(\eta_t^n)_t$, particle $k$ stops moving as soon as the particle is in position $n+1-k$ and particles $\{1, \ldots, k-1\}$ are to its right. Define the events (depending implicitly on $n$)
$$A_t^k = \{S(\nu_t^{k,n}, n-k) = k\}$$
$$= \{\text{particles } \{1, \ldots, k\} \text{ are in positions } \{n-k+1, \ldots, n\}\}.$$
The key to the result is the identity
$$(12) \qquad \{\beta^n(k) \leq t\} = A_t^k \cap A_t^{k-1}.$$

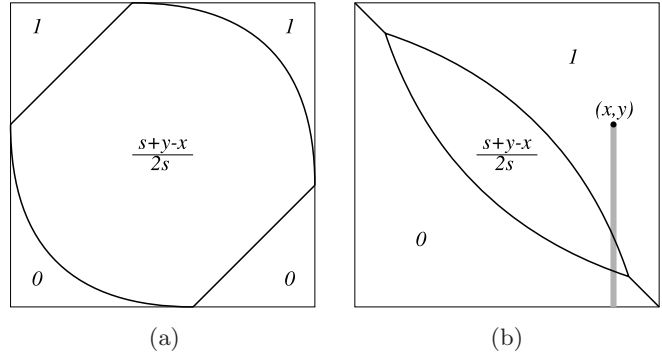

FIG. 8. *Values of $g(s,x,y)$ in the $(x,y)$ plane for* (a) $0 < s < 1$ *and* (b) $1 < s < 2$ *(there is no need to treat the boundary case $s = 1$ because of continuity in $s$). To compute $g(s,x,y)$, integrate $\kappa_s$ over the shaded thin strip of infinitesimal width.*

By Lemma 3.3, we can write

$$
\begin{aligned}
A_t^k &= \{S(B_n R_k \nu_t^k, n-k) = k\} \\
&= \{S(\nu_t^k, n-k) \wedge k \wedge (n-(n-k))^+ = k\} \\
&= \{S(\nu_t^k, n-k) \geq k\}.
\end{aligned}
\tag{13}
$$

Essentially, the fact that $S(\nu_t^k, \cdot)$ has a deterministic scaling limit implies that $A_t^k$ has a sharp threshold (that is, its probability increases from $\varepsilon$ to $1 - \varepsilon$ over an asymptotically small time interval). Thus, the finishing times are concentrated.

By symmetry, it suffices to prove that, with high probability, $\beta^n(k)$ is close to $\gamma_{k/n}$ for all $1 \leq k \leq n/2$. With this in mind, we start by showing this for a single particle. Let $k = k(n)$ be such that $k/n \to y \in [0, 1/2]$ as $n \to \infty$. From Theorem 4.3, we have, for any fixed $s$, that

$$
\begin{aligned}
\frac{1}{n} S(\nu_{ns}^k, n-k) &\xrightarrow[n\to\infty]{\mathbb{P}} \int_{1-y}^{\infty} h\left(\frac{x-y}{s}\right) dx \\
&= \begin{cases} 0, & 0 \leq s \leq 1-2y, \\ \dfrac{(s+2y-1)^2}{4s}, & s \geq 1-2y. \end{cases}
\end{aligned}
\tag{14}
$$

Thus, $A_{ns}^k$ and $A_{ns}^{k-1}$ hold with probability tending to 1 as soon as $\frac{(s+2y-1)^2}{4s} > y$, which simplifies to $s > 1 + 2\sqrt{y(1-y)} = \gamma_y$. Conversely, if $s < \gamma_y$ and $y \neq 0$, then $\mathbb{P}(A_{ns}^k)$ and $\mathbb{P}(A_{ns}^{k-1})$ tend to 0. This implies that $\beta^n(k)/n \xrightarrow{\mathbb{P}} \gamma_y$ for any $y \in (0, 1)$.

To get simultaneous convergence for all particles, we will use the easy facts that $S(\nu_{ns}^k, x)$ is increasing in $k$ and $s$, and decreasing in $x$. Since, for any fixed $s, y$, we have convergence in probability of $n^{-1} S(\nu_{ns}^k, n-k)$, this monotonicity implies that $n^{-1}S(\nu_{ns}^k, n-k)$ converges in probability to the right-hand side of (14), uniformly in $k$ and $s < 2$. [Proximity to the limit at a finite set of $(s,y)$ implies proximity at all intermediate points.]

What this implies for the finishing times is that for any $\varepsilon > 0$, the scaled finishing times $\beta^n(k)/n$ for particles $k \in \{\lfloor \varepsilon n \rfloor, \ldots, \lfloor n/2 \rfloor\}$ are (with asymptotically high probability) uniformly close to the given limit $\gamma_{k/n}$, that is,

$$\max_{\varepsilon n \leq k \leq n/2} |\beta^n(k)/n - \gamma_{k/n}| \xrightarrow[n\to\infty]{\mathbb{P}} 0.$$

For $k \in \{1, \ldots, \lfloor \varepsilon n \rfloor\}$, however, we only get an upper bound on the finishing times, namely, that

$$\max_{1 \leq k \leq \varepsilon n} (\beta^n(k)/n - \gamma_{k/n})^+ \xrightarrow[n\to\infty]{\mathbb{P}} 0.$$

Note that $\gamma_{k/n}$ is close to $\gamma_0 = 1$ for such small $k$, so we know that the particles with small labels must finish shortly after time $s = 1$, but not that they cannot finish much sooner. The reason we do not yet get a lower bound when $k \ll n$ is that for such $k$, the scaling limits of $k$ and $S(\nu_{n(1-\varepsilon)}^k, n-k)$ are both 0. Thus, the scaling limits are not enough to deduce that $S(\nu_{n(1-\varepsilon)}^k, n-k) < k$ and hence that the event $A_{ns}^k$ occurred. To complete the proof, it suffices to note that particle 1 performs a random walk, moving only to the right at random times with rate 1, and, therefore, at time $(1-2\varepsilon)n$, it is, with high probability, still at a distance of at least $\varepsilon n$ from position $n$. However, this implies that particles $\{1, \ldots, \lfloor \varepsilon n \rfloor\}$ are also not finished by time $(1-\varepsilon)n$. $\square$

PROOF OF THEOREM 1.6. For a configuration $\rho \in \Omega_0$, let $\pi_\rho(j)$ denote the position of the $j$th rightmost particle in $\rho$. By (12) and (13), we see that

$$\beta^n(k) = V^n(k) \vee V^n(k-1),$$

where

$$V^n(j) = \inf\{t > 0 : S(\nu_t^j, n-j) \geq j\}$$
$$= \inf\{t > 0 : \pi_{\nu_t^j}(j) = n+1-j\}.$$

If we were interested in $V^n(k)$, it would follow immediately from a theorem of Johansson [8], Theorem 1.6, that $V^n(k)$ converges in distribution after scaling to the Tracy–Widom distribution. Because our random time is a maximum of $V^n(k)$ and $V^n(k-1)$, we need to show that these two times cannot be very far apart. Define random times $W_1 = W_1^n$ and $W_2 = W_2^n$ by

$$W_1 = \inf\{t > 0 : \pi_{\nu_t^{k-1}}(k-1) = n+1-k\},$$
$$W_2 = \inf\{t > 0 : \pi_{\nu_t^{k-1}}(k) = n+1-k\}$$

and observe that we have the bounds
$$W_1 \leq \beta^n(k) = V^n(k) \vee V^n(k-1) \leq W_2.$$
The lower bound follows from the fact that
$$\pi_{\nu_t^k}(k) \leq \pi_{\nu_t^{k-1}}(k-1),$$
which implies that $V^n(k) \geq W_1$. Similarly, the upper bound follows from the inequalities
$$\pi_{\nu_t^{k-1}}(k) < \pi_{\nu_t^{k-1}}(k-1) \quad \text{and} \quad \pi_{\nu_t^{k-1}}(k) \leq \pi_{\nu_t^k}(k+1) < \pi_{\nu_t^k}(k),$$
which imply, respectively, that $V^n(k-1) \leq W_2$ and $V^n(k) \leq W_2$. Now, use Johansson's theorem ([8], Theorem 1.6, see also [8], Corollary 1.7) for $W_1$ and $W_2$, to get that

(15) $$\frac{W_i - \gamma_{k/n} n}{\gamma_y^{2/3}(y(1-y))^{-1/6} n^{1/3}} \Longrightarrow F_{\mathrm{TW}} \quad \text{as } n \to \infty \ (i = 1, 2)$$

(we use the well-known equivalence between the TASEP with Rost's step initial conditions $\mathbf{1}_{(-\infty, 0]}$ and directed last-passage percolation in $\mathbb{N}^2$; in the notation of the paper [8], we have $N = yn$ and $\gamma = (1-y)/y$). Since $\beta^n(k)$ is bounded between two random variables having the same distributional limit (with the same scaling), it must also converge in distribution to $F_{\mathrm{TW}}$ with the same scaling. $\square$

**6. The inversion number.** In this section, we prove Theorem 1.4. By Theorem 1.5, it is enough to consider $s \leq 2$ since, for $s > 2$, we have, with high probability, that $\eta_{ns}^n$ is the reverse permutation and, in particular, $\mathrm{inv}(\eta_{ns}^n) = \binom{n}{2}$.

Let $\eta \in \mathcal{S}_n$ be any permutation with normalized empirical measure $\mu$ defined as in (1). Thus, $\mu$ is supported within the square $[0,1]^2$. For points $z = (x, y)$ and $z' = (x', y')$ in $[0, 1]^2$, we write $\{z \searrow z'\}$ if $\{x < x' \text{ and } y > y'\}$. The basic observation is that we can express the inversion number $\mathrm{inv}(\eta)$ in terms of the measure $\mu$, as follows. If we sample independent random points $z, z' \in [0, 1]^2$ with distribution $\mu$, then
$$\frac{1}{n^2} \mathrm{inv}(\eta) = \mu \otimes \mu(z \searrow z') = \int_{[0,1]^2} \int_{[0,1]^2} \mathbf{1}[z \searrow z'] \mu(dz') \mu(dz).$$
The integral above, as a function of the measure $\mu$, is continuous at all measures that assign 0 measure to all horizontal and vertical lines (since, for such $\mu$, the product measure $\mu \otimes \mu$ will assign 0 measure to the boundary of the set $\{z \searrow z'\} \subset [0, 1]^4$). This includes, in particular, the measures $\kappa_s$. By Theorem 1.2, it follows that
$$\binom{n}{2}^{-1} \mathrm{inv}(\eta_{ns}^n) \xrightarrow[n \to \infty]{\mathbb{P}} I(s) := \int_{[0,1]^2} \int_{[0,1]^2} 2 \cdot \mathbf{1}[z \searrow z'] \kappa_s(dz') \kappa_s(dz).$$

Thus, calculating the limit of the scaled number of inversions function is reduced to evaluating the $\kappa_s \otimes \kappa_s$ measure of a certain set in $[0,1]^4$. This integration is tricky, so we sketch an argument below.

We first give a geometric description of $\kappa_s$ which will shed some light on subsequent formulae. To sample a point $(x,y)$ from $\kappa_s$, first choose $y$ uniformly in $[0,1]$ and $x$ uniformly in $[y-s, y+s]$ [so that $(x,y)$ is uniform in a parallelogram of area $2s$]. Next, evaluate $L_y^-(s)$ and $L_y^+(s)$. If $x$ is outside the interval $[L_y^-, L_y^+]$, replace it by whichever of $L_y^-, L_y^+$ is closer to $x$. When $s > 1$, if $L^- > L^+$, then replace $x$ by $1 - y$ (see Figure 9). The equivalence of this description to the original definition of $\kappa_s$ can be verified with a simple computation that we omit (it also has a more conceptual explanation in terms of Theorems 1.2 and 4.3 and the combinatorial operations $B_n$ and $R_k$; e.g., look again at Figure 7).

CLAIM 6.1. *For a point $z = (x,y)$ in the support of $\kappa_s$, if $z'$ has law $\kappa_s$, then*

$$\kappa_s(z \searrow z') = \begin{cases} (y - x + s)^2/4s, & s \leq 1 \text{ or } y \in [W_s^-, W_s^+], \\ y, & s > 1 \text{ and } y \notin [W_s^-, W_s^+]. \end{cases}$$

PROOF. The second-case is trivial. The first case is equivalent to showing that the part of the parallelogram with $z'$ satisfying $z \searrow z'$ has area $(y - x + s)^2/2$. This region is a fairly simple polygon [see Figure 10(a)] and its area is easily computed.

When $s \leq 1$, there is a nice geometric argument. Since $\kappa_s$ is symmetric with respect to the reflection $(x,y) \to (y,x)$ (check directly, or deduce it from Theorems 1.2 and 1.3), the symmetric description of how to sample $z'$ by choosing $x'$ first and then $y'$, and shifting $y'$ vertically if the resulting point is outside the support of $\kappa_s$ is equally valid. Therefore, we can reflect the curved triangle $ABF$ of Figure 9 and deduce that $\kappa_s(z' : z \searrow z')$, which is equal to the area of the shaded region in Figure 10(a), is also the area of the shaded region in Figure 10(b), which is simply a right triangle.

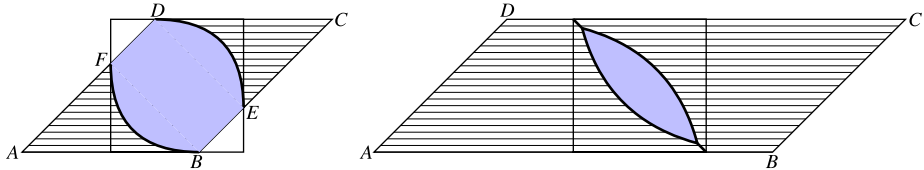

FIG. 9. *The measure $\kappa_s$: a point is selected uniformly in the parallelogram $ABCD$. If it is outside the curved polygon $BEDF$ where $\kappa_s$ is supported, then it is shifted horizontally into the support. A similar description also applies for $s > 1$.*

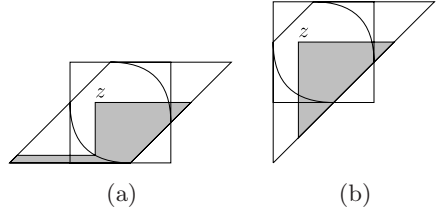

FIG. 10. *Proof of the claim when $s \leq 1$. We will have $z' \searrow z$ if the point from which $z'$ is obtained by shifting the x-coordinate lies in the shaded region in* (a). *Symmetry of $\kappa_s$ implies that the areas of the shaded regions in* (a) *and* (b) *are equal.*

For the case $s > 1$, we may either deduce the claim from the case $s < 1$ together with analyticity of the areas involved, or prove it directly by computing the area of the hexagon, as in Figure 10(a). □

We now compute $I(s)$, starting with the case $s \leq 1$. From the claim, it follows that, in this case,

$$I(s) = \int_{[0,1]^2} 2\frac{(y-x+s)^2}{4s} \kappa_s(dz).$$

Here, the contribution from $z$ on the lower arc [with $x = L_y^-(s)$] is

$$\int_0^s \int_{y-s}^{L_y^-(s)} 2\frac{(y-L_y^-(s)+s)^2}{4s} \frac{dx}{2s}\, dy = \frac{s^2}{5}.$$

The upper arc gives an equal contribution.

To integrate $(y-x+s)^2/4s$ in the interior, note that, for any $s \leq 1/2$, the integral of a function $G$ in the internal region is given by

(16) $$\int_0^s \int_{L^-}^{y+s} G\,dx\,dy + \int_s^{1-s} \int_{y-s}^{y+s} G\,dx\,dy + \int_{1-s}^1 \int_{y-s}^{L^+} G\,dx\,dy.$$

For $1/2 < s \leq 1$, it seems, at first glance, that a different formula is required, but, in fact, it is easy to see that (16) is still correct due to cancellation. Applying this with $G = (y-x+s)^2/(8s^2)$ and adding the arc contributions gives the answer $I(s) = \frac{2}{3}s - \frac{1}{15}s^2$ for $s \leq 1$.

The case $s > 1$ is trickier. The contributions from $z$ with $y < W_s^-$ or $y > W_s^+$ are simple enough. The integrals for the arcs and interior are also not difficult to evaluate, but, due to having $W_s^\pm$ as endpoints, give less elegant expressions. Overall, the different parts shown in Figure 11 evaluate to

$$A + A' = 1 - X,$$
$$B = B' = X + \frac{1}{5\sqrt{2}} s^{-1/2}((1-X)^{5/2} - (1+X)^{5/2}),$$

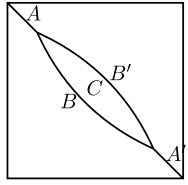

FIG. 11. *Contribution to $I(s)$ for $s > 1$. $A, A'$: diagonal segments; $B, B'$: arcs; $C$: lens-shaped region between the arcs.*

$$C = \left(-1 - \frac{2}{3}s\right)X - \frac{1}{15\sqrt{2}}s^{-1/2}((1-X)^{5/2} - (1+X)^{5/2})$$
$$- \frac{1}{15\sqrt{2}}[s^{-1/2}(1 + 12s - s^2)Y + (2+10s)(2-s)^{1/2}Z],$$

where we denote

$$X = \sqrt{s(2-s)}, \qquad Y = \sqrt{1-X} - \sqrt{1+X},$$
$$Z = \sqrt{1-X} + \sqrt{1+X}.$$

It seems impossible to simplify this any further until one realizes that $Y = -\sqrt{2(2-s)}$ and $Z = \sqrt{2s}$. Then, straightforward manipulation of $I(s) = A + A' + B + B' + C$ gives the theorem.

**7. Particle trajectories and second-class particles.** Our goal in this section is to prove Theorem 1.1. To this end, note that one can identify the location of particle $k$ in the oriented swap process $(\eta_t^n)_t$ at a given time $t$ as the unique place where the two configurations $\nu_t^{k,n}$ and $\nu_t^{k-1,n}$ differ. Because of Lemma 3.3, this is the same as the place where the two configurations $B_n R_k \nu_t^k$ and $B_n R_{k-1} \nu_t^{k-1}$ differ.

Note that $\nu_t^k, \nu_t^{k-1}$ also differ in a unique place. This is true at time $t = 0$, and it is easy to see that it remains true whenever one of the jump operators $J_m$ is applied. In the context of the TASEP, the place where $\nu_t^k, \nu_t^{k-1}$ differ is referred to as a *second-class particle* since it behaves in effect like a particle which has priority over holes (i.e., can swap with a hole to its right and does so with exponential rate 1), but over which the "first-class" particles initially at $(-\infty, k-1]$ have priority (to see this, try to imagine the effect that the operator $J_m$ has on the pair of configurations).

As before, let $\Omega_0$ denote the space of configurations in $\{0,1\}^{\mathbb{Z}}$ having a rightmost particle and, if a configuration is in $\{0,1\}^I$ for some subinterval $I \subset \mathbb{Z}$, consider it to be in $\{0,1\}^{\mathbb{Z}}$ by extending its definition to be 0 outside $I$. We call two configurations $\rho, \rho' \in \Omega_0$ *compatible* if $\rho$ and $\rho'$ differ in exactly one position where $0 = \rho'(x) < \rho(x) = 1$, and we denote the position where they differ by $\Sigma_{\rho,\rho'}$. For $\rho \in \Omega_0$, denote the position of the $j$th rightmost particle of $\rho$ by $\pi_\rho(j)$ and denote by $\theta_\rho(n,j)$ the position of the $j$th rightmost hole of $\rho \vee \mathbf{1}_{[n+1,\infty)}$ (in words, the $j$th rightmost hole among the holes that are in $(-\infty, n]$; define it as $n+1$ if $j=0$). The following lemma elucidates the effect that the transformations $B_n$ and $R_k$ have on the second-class particle.

LEMMA 7.1. *Let $\rho, \rho' \in \Omega_0$ be compatible configurations. Then the pair $R_k\rho, R_{k-1}\rho'$ is compatible, with its second-class particle in position*

$$(17) \qquad \Sigma_{R_k\rho, R_{k-1}\rho'} = \Sigma_{\rho, \rho'} \vee \pi_\rho(k).$$

*If $\rho$ has infinitely many holes in $(-\infty, 0]$, then the pair $B_n\rho, B_n\rho'$ is also compatible, with a second-class particle in*

$$(18) \qquad \Sigma_{B_n\rho, B_n\rho'} = \Sigma_{\rho, \rho'} \wedge \theta_\rho(n, S(\rho, n)).$$

*Consequently, for $k \leq n$ and any compatible $\rho, \rho' \in \Omega_0$, the pair of configurations $B_n R_k \rho, B_n R_{k-1}\rho'$ is compatible with a second-class particle in*

$$(19) \qquad \Sigma_{B_n R_k \rho, B_n R_{k-1}\rho'} = \Sigma_{\rho, \rho'} \vee \pi_\rho(k) \wedge \theta_{R_k\rho}(n, S(R_k\rho, n)).$$

PROOF. Note that (19) follows from (17) and (18) [since $R_k\rho$ satisfies the assumption of (18)]. For (17), note that if $\Sigma_{\rho, \rho'} > \pi_\rho(k)$, then, when applying $R_k$ to $\rho$ and also when applying $R_{k-1}$ to $\rho'$, all particles to the left of $\pi_\rho(k) = \pi_{\rho'}(k-1)$ will be deleted, so the transformed configurations still differ at $\Sigma_{\rho, \rho'}$. If $\Sigma_{\rho, \rho'} \leq \pi_\rho(k)$, then, in $R_k\rho$, all particles to the left of $\pi_\rho(k)$ are deleted and in $R_{k-1}\rho'$, all particles to the left of $\pi_{\rho'}(k-1) = \pi_\rho(k-1)$ are deleted, so the two transformed configurations now differ at $\pi_\rho(k)$. This proves (17).

To prove (18), first note that an equivalent description of the transformation $B_n$ acting on $\rho$ is that it takes the $S(\rho, n)$ particles to the right of $n$ and places them in the $S(\rho, n)$ rightmost holes to the left of $n+1$. This is true for configurations with sufficiently many holes to the left of $n+1$, which is why, for convenience, we assumed that there are infinitely many such holes in $\rho$. Similarly, $B_n$ acting on $\rho'$ has the effect of filling in the $S(\rho', n)$ rightmost holes in $\rho'$ to the left of $n+1$. Now, assume that $\Sigma_{\rho, \rho'} < \theta_\rho(n, S(\rho, n))$. In particular, it follows that $\Sigma_{\rho, \rho'} < n+1$, so $S(\rho, n) = S(\rho', n)$. Because of the assumption, the $S(\rho, n)$ rightmost holes to the left of $n+1$ are the same for $\rho$ and $\rho'$, and they both get filled under the operation $B_n$. So, $B_n\rho$ and $B_n\rho'$ still differ only at position $\Sigma_{\rho, \rho'}$ and (18) is correct. The other possibility is that $\Sigma_{\rho, \rho'} > \theta_\rho(n, S(\rho, n))$ [we cannot have equality since $\rho$ has a particle at position $\Sigma_{\rho, \rho'}$ and a hole at position $\theta_\rho(n, S(\rho, n))$]. For this case, consider two further sub-cases: first, if $\Sigma_{\rho, \rho'} \leq n$, then $S(\rho, n) = S(\rho', n)$, and then the holes in $\rho$ in positions $\theta_\rho(n, 1), \theta_\rho(n, 2), \ldots, \theta_\rho(n, S(\rho, n))$ get filled in $B_n\rho$, whereas for $\rho'$, the holes in positions $\theta_\rho(n, 1), \ldots, \theta_\rho(n, S(n, \rho) - 1), \Sigma_{\rho, \rho'}$ get filled when applying $B_n$. It follows that $B_n\rho$ and $B_n\rho'$ differ in position $\theta_\rho(n, S(\rho, n))$, as claimed. Finally, if $\Sigma_{\rho,\rho'} > n$, then $S(\rho', n) = S(\rho, n) - 1$, and then, in $\rho$, the holes in positions $\theta_\rho(n, 1), \ldots, \theta_\rho(n, S(\rho, n))$ get filled, but in $\rho'$, the holes in positions $\theta_\rho(n, 1), \ldots, \theta_\rho(n, S(\rho, n) - 1)$ get filled, so, again, $B_n\rho$ and $B_n\rho'$ differ in position $\theta_\rho(n, S(n, \rho))$. $\square$

Denote by $X_t$ the location of the second-class particle in the TASEP pair $\nu_t^0, \nu_t^{-1}$. Guiol and Mountford [11] proved that $X_t/t$ converges a.s. as $t \to \infty$ to a random variable distributed uniformly in $[-1, 1]$. The following result is an easy consequence of their result. A similar, but slightly weaker, result was proven earlier by Ferrari and Kipnis [7]: there, the limit is with respect to the (weaker) topology of pointwise convergence of functions.

THEOREM 7.2 (Guiol and Mountford [11]). *Let $U$ denote a random variable distributed uniformly on $[-1, 1]$. Let $\hat{X}^n(s) = n^{-1} X_{ns}$ denote the trajectory of the second-class particle when space and time are scaled by $n$, considered as a random function in the function space $\mathbb{R}^{[0,\infty)}$. Then,*

$$\hat{X}^n(s) \Longrightarrow U \cdot s.$$

*Here, "$\Longrightarrow$" denotes convergence in distribution with respect to the uniform topology on functions $[0, \infty) \to \mathbb{R}$.*

PROOF OF THEOREM 1.1. Let $k = k(n)$ be such that $k/n \to y \in [0, 1]$ when $n \to \infty$. For $s > \gamma_y$, we already know, from Theorem 1.5, that with asymptotically high probability, $T_k^n(s) = \frac{n+1-k}{n} \to 1 - y = \phi_y(s)$, so it will be enough to prove the claimed convergence on the space of functions $\mathbb{R}^{[0,\gamma_y]}$. From Lemmas 3.3 and 7.1, we have that

$$T_k^n(s) = \frac{1}{n} \Sigma_{\nu_{ns}^{k,n}, \nu_{ns}^{k-1,n}} = \frac{1}{n} \Sigma_{B_n R_k \nu_{ns}^k, B_n R_{k-1} \nu_{ns}^{k-1}}$$
$$= \frac{1}{n}[\Sigma_{\nu_{ns}^k, \nu_{ns}^{k-1}} \vee \pi_{\nu_{ns}^k}(k) \wedge \theta_{R_k \nu_{ns}^k}(n, S(\nu_{ns}^k, n) \wedge k)].$$

Therefore, we need to understand each of the components on the right-hand side. Note that, by translation, $\Sigma_{\nu_{ns}^k, \nu_{ns}^{k-1}}$ is equal in distribution to $k + \Sigma_{\nu_{ns}^0, \nu_{ns}^{-1}} = k + \hat{X}^n(s)$. Therefore, by Theorem 7.2, we know that

$$n^{-1} \Sigma_{\nu_{ns}^k, \nu_{ns}^{k-1}} \Longrightarrow y + Us.$$

We now claim that

(20) $$n^{-1} \pi_{\nu_{ns}^k}(k) \xrightarrow[n \to \infty]{\mathbb{P}} \Lambda_y^-(s),$$

$$n^{-1} S(\nu_{ns}^k, n) \xrightarrow[n \to \infty]{\mathbb{P}} \psi_y(s) := \int_1^\infty h\left(\frac{x-y}{s}\right) dx$$

(21)
$$\stackrel{\text{(see Fig. 7)}}{=} \begin{cases} 0, & 0 \leq s < 1-y, \\ \dfrac{(s+y-1)^2}{4s}, & s \geq 1-y, \end{cases}$$

(22) $\theta_{R_k \nu_{ns}^k}(n, S(\nu_{ns}^k, n) \wedge k) \xrightarrow[n \to \infty]{\mathbb{P}} \Lambda_y^+(s).$

If we prove these claims, it will follow that

(23) $\qquad \hat{T}_k^n \Longrightarrow (y + Us) \vee \Lambda_y^-(s) \wedge \Lambda_y^+(s) = \phi_y(s) \qquad (s \leq \gamma_y),$

which is the claim of the theorem. Relation (21) follows immediately from Corollary 4.2. The other two relations, (20) and (22), are also relatively straightforward consequences of Corollary 4.2. To prove (20), note that $\pi_{\nu_{ns}^k}(k)$ satisfies

$$S(\nu_{ns}^k, \pi_{\nu_{ns}^k}(k)) = k - 1$$

[indeed, it is the minimal $x$ with $S(\nu_{ns}^k, x) = k - 1$]. On the other hand, if $y > 0$, then, by Corollary 4.2, for any $\varepsilon > 0$, we have that, as $n \to \infty$,

(24)
$$n^{-1} S(\nu_{ns}^k, (\Lambda_y^-(s) + \varepsilon)n)$$
$$\xrightarrow{\mathbb{P}} \int_{\Lambda_y^-(s)+\varepsilon}^{\infty} h\left(\frac{x-y}{s}\right) dx < y - \delta \quad \text{and}$$

(25)
$$n^{-1} S(\nu_{ns}^k, (\Lambda_y^-(s) - \varepsilon)n)$$
$$\xrightarrow{\mathbb{P}} \int_{\Lambda_y^-(s)-\varepsilon}^{\infty} h\left(\frac{x-y}{s}\right) dx > y + \delta$$

for some $\delta = \delta(y, \varepsilon) > 0$ that depends on $\varepsilon$. By monotonicity of $S(\rho, x)$ in $x$, it follows that

$$\mathbb{P}[\Lambda_y^-(s) - \varepsilon < n^{-1} \pi_{\nu_{ns}^k}(k) < \Lambda_y^-(s) + \varepsilon] \to 1 \qquad \text{as } n \to \infty,$$

which implies (20) since $\varepsilon$ was arbitrary.

In the extremal case $y = 0$ (which, by symmetry, also implies the case $y = 1$), (25) still holds, giving the lower bound for $\pi_{\nu_{ns}^k}(k)$, and (24) does not, but we instead use the easy fact that

$$n^{-1} \pi_{\nu_{ns}^k}(k) < n^{-1} \pi_{\nu_{ns}^k}(1) \xrightarrow{\mathbb{P}} s = \Lambda_0^-(s) \qquad \text{as } n \to \infty$$

(since the rightmost particle performs a random walk jumping to the right at rate 1) to get the upper bound, so (20) still holds.

Next, to prove (22), we first simplify the left-hand side by noting the easily checked fact that $\psi_y(s) \leq y$ for $0 \leq s \leq 1 + y + 2\sqrt{y}$ and, in particular, for $0 \leq s \leq \gamma_y$ since $\gamma_y \leq 1 + y + 2\sqrt{y}$. Therefore, we have

$$\frac{1}{n}(S(\nu_{ns}^k, n) \wedge k) \xrightarrow[n \to \infty]{\mathbb{P}} \psi_y(s).$$

Define

$$Z_{y,s}(x) = \begin{cases} 0, & x < \Lambda_y^-(s), \\ h\left(\dfrac{x-y}{s}\right), & x > \Lambda_y^-(s). \end{cases}$$

By Corollary 4.2 and (20), $Z_{y,s}(x)$ is the limiting density profile of the process $R_k \nu_{ns}^k$, in the sense that, for all $u \in \mathbb{R}$, we have

$$n^{-1} S(R_k \nu_{ns}^k, j_n) \xrightarrow[\substack{n \to \infty \\ j_n/n \to u}]{\mathbb{P}} \int_u^\infty Z_{y,s}(x)\, dx.$$

It follows, by an argument similar to the one used to prove (20) above, that

$$n^{-1} \theta_{R_k \nu_{ns}^k}(n, S(\nu_{ns}^k, n) \wedge k)$$
$$\xrightarrow[n \to \infty]{\mathbb{P}} \inf\left\{ u : \int_u^1 (1 - Z_{y,s}(x))\, dx \geq \psi_y(s) \right\}.$$

For $s \leq \gamma_y$, this last function is easily seen (refer again to Figure 7) to equal $\Lambda_y^+(s)$. This completes the proof of (22) and therefore also of (23). □

## 8. Additional comments and open problems.

1. *Uniformly random sorting networks.* The uniform sorting network model of [2] exhibits behavior similar to the oriented swap process. A key problem is to prove the conjectures in [2] that are the analogs of Theorems 1.1 and 1.2. One possible approach would be to try to relate uniform sorting networks to some variant of the random swap process which can be analyzed using the exclusion process techniques developed in this paper.

2. *Limiting distribution of the absorbing time.* Theorem 1.6 gives the limiting distribution of the fluctuations of the finishing times of individual particles. However, the relation between finishing times of different particles is more delicate and requires knowledge about the joint distribution of last-passage percolation times. An interesting open problem would be to find sequences of scaling constants $(a_n)_{n=1}^\infty$, $(b_n)_{n=1}^\infty$ and a distribution function $F$ such that the absorbing time $\beta_*^n$ of the oriented swap process satisfies the convergence in distribution

$$a_n(\beta_*^n - 2n) - b_n \Longrightarrow F \qquad \text{as } n \to \infty.$$

3. *Partially asymmetric swap processes.* The *asymmetric exclusion process (ASEP)* is defined similarly to the TASEP, with the difference that particles can jump in either direction, where jumps to the right happen at rate $\alpha$ and to the left at rate $\alpha'$. When $\alpha' = 0$, this is (a time-change of) the TASEP. When $\alpha > \alpha' > 0$, this is the partially asymmetric exclusion process (PASEP).

Some results for the TASEP have known analogs for the PASEP, while others do not. In particular, the following conjecture is needed to prove an analogue of Theorem 1.1 for the partially asymmetric swap process mentioned in the Introduction:

CONJECTURE 8.1. Theorem 7.2 holds for the partially asymmetric model with parameters $\alpha > \alpha'$, with speed $U$ uniform on $[\alpha' - \alpha, \alpha - \alpha']$.

O. ANGEL
DEPARTMENT OF MATHEMATICS
UNIVERSITY OF BRITISH COLUMBIA
121-1984 MATHEMATICS ROAD
VANCOUVER, BRITISH COLUMBIA
CANADA V6T 1Z2
E-MAIL: angel@math.utoronto.ca

A. E. HOLROYD
MICROSOFT RESEARCH
1 MICROSOFT WAY
REDMOND, WASHINGTON 98052
USA
AND
DEPARTMENT OF MATHEMATICS
UNIVERSITY OF BRITISH COLUMBIA
121-1984 MATHEMATICS ROAD
VANCOUVER, BRITISH COLUMBIA
CANADA V6T 1Z2
E-MAIL: holroyd@math.ubc.ca

D. ROMIK
EINSTEIN INSTITUTE OF MATHEMATICS
THE HEBREW UNIVERSITY
GIVAT-RAM
JERUSALEM 91904
ISRAEL
E-MAIL: romik@math.huji.ac.il